\documentclass[a4paper,12pt]{article}

\usepackage{amscd}
\usepackage{amssymb,amsfonts,amsmath,amsthm}
\usepackage{verbatim}

\usepackage{amsmath}
\usepackage{amsthm}
\usepackage{amssymb}

\usepackage{latexsym}
\usepackage{array}
\usepackage{enumerate}
\usepackage{bm}

\numberwithin{equation}{section}

\usepackage{amsmath,amssymb,mathrsfs,amsthm}
\usepackage[all]{xy}
\usepackage{graphicx}
\usepackage{ulem}
\usepackage{ascmac}
\usepackage{color}
\usepackage{mathtools}

\newcommand{\rHom}{{\bfR}\mathrm{Hom}}
\newcommand{\BDC}{{\mathbf{D}}^{\mathrm{b}}}

\newcommand{\Mod}{\mathrm{Mod}}
\newcommand{\Hom}{\mathrm{Hom}}

\newcommand{\CC}{\mathbb{C}}
\newcommand{\NN}{\mathbb{N}}
\newcommand{\RR}{\mathbb{R}}

\newcommand{\ZZ}{\mathbb{Z}}
\newcommand{\D}{\mathcal{D}}

\newcommand{\M}{\mathcal{M}}

\newcommand{\calT}{\mathcal{T}}

\newcommand{\calS}{\mathcal{S}}

\newcommand{\supp}{{\rm supp}}
\renewcommand{\SS}{{\rm SS}}

\newcommand{\Int}{{\rm Int}}

\newcommand{\Perv}{{\rm Perv}}
\newcommand{\Sol}{{\rm Sol}}

\newcommand{\simto}{\overset{\sim}{\longrightarrow}}
\newcommand{\simot}{\overset{\sim}{\longleftarrow}}

\newcommand{\op}{{\mbox{\scriptsize op}}}

\newcommand{\SO}{\mathcal{O}}

\newcommand{\SM}{\mathcal{M}}

\newcommand{\Modhol}{\mathrm{Mod}_{\mbox{\rm \scriptsize hol}}}
\newcommand{\Modrh}{\mathrm{Mod}_{\mbox{\rm \scriptsize rh}}}

\newcommand{\BDCcoh}{{\mathbf{D}}^{\mathrm{b}}_{\mbox{\rm \scriptsize coh}}}

\newcommand{\BDChol}{{\mathbf{D}}^{\mathrm{b}}_{\mbox{\rm \scriptsize hol}}}
\newcommand{\BDCrh}{{\mathbf{D}}^{\mathrm{b}}_{\mbox{\rm \scriptsize rh}}}

\newcommand{\rhom}{{\bfR}{\mathcal{H}}om}

\newcommand{\I}{{\rm I}}

\newcommand{\var}[1]{\overline{#1}}

\newcommand{\bfR}{\mathbf{R}}

\newcommand{\bfD}{\mathbf{D}}

\newcommand{\rmt}{{\rm t}}

\newcommand{\reg}{{\rm reg}}

\newcommand{\CCc}{\CC\mbox{\scriptsize -}c}
\newcommand{\RRc}{\RR\mbox{\scriptsize -}c}
\newcommand{\frCCc}{{\rm fr}\mbox{\scriptsize -}\CC\mbox{\scriptsize -}c}
\newcommand{\frRRc}{{\rm fr}\mbox{\scriptsize -}\RR\mbox{\scriptsize -}c}
\newcommand{\wRRc}{{\rm w}\mbox{\scriptsize -}\RR\mbox{\scriptsize -}c}
\newcommand{\wCCc}{{\rm w}\mbox{\scriptsize -}\CC\mbox{\scriptsize -}c}
\newcommand{\wc}{{\rm w}\text{\scriptsize{-}}\calT\text{\scriptsize{-}c}}

\newcommand{\Char}{\rm Char}

\newtheorem{theorem}{Theorem}[section]
\newtheorem{conjecture}[theorem]{Conjecture}
\newtheorem{corollary}[theorem]{Corollary}
\newtheorem{lemma}[theorem]{Lemma}

\newtheorem{proposition}[theorem]{Proposition}

\newtheorem*{proposition*}{Proposition}

\theoremstyle{definition}
\newtheorem{definition}[theorem]{Definition}
\theoremstyle{remark}

\newtheorem{maintheorem}{{\bf Theorem}}

\title{{\Large Regularity Criteria for Ind-Sheaves}}
\author{
Yohei Ito\thanks{
Department of Mathematical Sciences,
Aoyama Gakuin University, 
5-10-1 Fuchinobe, 
Chuo-ku, 
Sagamihara-shi, 
Kanagawa,
252-5258, 
Japan.
E-mail: yitoh@math.aoyama.ac.jp}
\and 
Ryosuke Sakamoto\thanks{
Institute for the Advanced Study of Human Biology (ASHBi),
Kyoto University,
Yoshida-Konoe-cho,
Sakyo-ku,
Kyoto,
606-8501,
Japan.
E-mail: sakamoto.ryosuke.7x@kyoto-u.ac.jp} 
}

\date{}

\begin{document}
\maketitle

\begin{abstract}
In this paper, we prove the Kashiwara--Schapira regularity conjecture \cite[Conj.\:6.4]{KS03} in arbitrary dimension.
For this purpose, we show that under suitable finiteness and isotropic micro-support conditions, ind-sheaves which are regular on the zero section come from $\RR$-constructible sheaves.
Furthermore, we characterize the essential image of constructible sheaves under the natural embedding functor, and relate them to regular holonomic $\D$-modules via the regular Riemann--Hilbert correspondence.
Moreover, we present a non-characteristic deformation lemma and a microlocal Morse lemma for regular ind-sheaves as an application of the main theorem.
\end{abstract}

%\noindent
%\textbf{Keywords}: irregular Riemann--Hilbert correspondence, equivariant sheaf, enhanced subanalytic sheaf.\\
%\textbf{MSC2020}: primary 32C38; secondary, 18F10, 35Q15, 32S60.

\section{Introduction}
Throughout this paper, 
$\Bbbk$ is a field and $X$ is a real analytic or complex manifold.

In \cite{KS01}, M.\:Kashiwara and P.\:Schapira introduced the notion of ind-sheaves
to treat ``sheaves" of functions with tempered growth conditions.
We denote by $\I\Bbbk_X$ the category of ind-sheaves of $\Bbbk$-vector spaces on $X$
and by $\beta_X\dashv\alpha_X\dashv\iota_X$
the standard functors between the category $\Mod(\Bbbk_X)$ of sheaves of $\Bbbk$-vector spaces on $X$
and $\I\Bbbk_X$ in \cite{KS01}.
In particular, we obtain the canonical morphism
$$u_F\colon F\to \iota_X\alpha_X(F)$$
for $F\in \BDC(\I\Bbbk_X)$.
Here, $\BDC(\I\Bbbk_X)$ is the derived category of $\I\Bbbk_X$.
Moreover, the micro-support of ind-sheaves was introduced in \cite[Def.\:4.2]{KS03}, and regularity for ind-sheaves was introduced in \cite[Def.\:4.9]{KS03} (cf. Definition \ref{regularity}).

In this framework,
the Kashiwara--Schapira regularity conjecture was proposed in \cite[Conj.\:6.4]{KS03}.
Let $X$ be a complex manifold.
We denote by $\Modhol(\D_X)$ (resp.~$\Modrh(\D_X)$) the abelian category of
holonomic (resp.~regular holonomic) $\D_X$-modules and by $\Sol_X$ (resp.~$\Sol^{\rm t}_X$)
the solution functor (resp.~the tempered solution functor) on $X$.
We refer to \cite{Bjo93, HTT} for the theory of $\D$-modules and to \cite[\S\:7]{KS01}
for the tempered solution functor.

\begin{conjecture}[{\cite[Conj.\:6.4]{KS03}}]
For any $\SM\in \Modhol(\D_X)$,
$\M\in \Modrh(\D_X)$ if and only if $\Sol^{\rm t}_X(\M)$ is regular.
\end{conjecture}

By \cite[Lem.\:7.4.11]{KS01}, for any $\M\in \BDCrh(\D_X)$, the canonical morphism
$$u_{\Sol_X^{\rmt}(\M)}\colon \Sol^{\rm t}_X(\M)\to \iota_X\Sol_X(\M)$$
is an isomorphism.
Therefore, $\Sol^{\rm t}_X(\M)$ is regular by \cite[Prop.\:4.10 (iii)]{KS03}.
%\cite[Prop.\:4.10 (iii)]{KS03}.

The converse was proved by \cite{Mar} in dimension one
by reducing to irregular holonomic $\D$-modules of the form $\D^m/\D^m P$,
showing that an irregular singularity is an obstruction to regularity
at the corresponding point on the zero section.

In this paper, we prove the converse in arbitrary dimension and in the derived setting.
More precisely, we prove the following main theorem.
We denote by $\BDChol(\D_X)$ (resp.~$\BDCrh(\D_X)$) the triangulated category of
holonomic (resp.~regular holonomic) $\D_X$-modules.
\begin{maintheorem}\label{mainA}
Let $\M\in \BDChol(\D_X)$.
The following are equivalent:
\begin{itemize}
\item[(i)]
$\M\in \BDCrh(\D_X)$.

\item[(ii)]
$\Sol^{\rm t}(\M)$ is regular.

\item[(iii)]
For any $x\in X$, $\Sol^t(\M)$ is regular along $\Char(\M)$ at $(x;0)\in T^\ast X$.

\item[(iv)]
The canonical morphism $\Sol^{\rm t}_X(\M)\to \iota_X\Sol_X(\M)$ is an isomorphism.
\end{itemize}
\end{maintheorem}

As already seen, (i) implies (ii), and it is obvious that (ii) implies (iii).
We prove that (iv) implies (i) by using the regular Riemann--Hilbert correspondence,
and (iii) implies (iv) via the following main theorem.
\begin{maintheorem}\label{mainB}
Let $X$ be a real analytic manifold,
$U_0\subset X$ be an open neighborhood of $x\in X$
and $\Lambda\subset T^\ast U_0$ be a conic closed subanalytic isotropic subset.
Suppose that $F\in \BDC(\I\Bbbk_X)$ is regular along $\Lambda$ at $(x; 0)$ and that
$$F|_{U_0}\in \BDC_{\I\RRc}(\I\Bbbk_{U_0}),\qquad
 \dim_k H^q(\alpha_XF)_y<\infty
 \quad(\forall y\in U_0,\ \forall q\in\mathbb Z).$$
Then, there is an open neighborhood $U_x\subset U_0$ of $x$ such that
$$u_F|_{U_x}\colon F|_{U_x}\xrightarrow{\ \sim\ }(\iota_X\alpha_XF)|_{U_x},\qquad
(\alpha_XF)|_{U_x}\in \BDC_{\RRc,\Lambda}(\Bbbk_{U_x}).$$
Here, the subscript $\Lambda$ means micro-support contained in $\Lambda\cap T^\ast U_x$.
\end{maintheorem}

Thus, our approach, using this theorem, together with the adjoint pair $\rhom(\cdot, \SO_X^\rmt) \dashv \Sol^{\rm t}_X$,
proves the conjecture in arbitrary complex dimension
without reducing to irregular holonomic $\D$-modules of the form $\D^m/\D^m P$.

We also describe the essential images of $\BDC_{\RRc}(\Bbbk_X)$ and $\BDC_{\CCc}(\CC_X)$ under the functor $\iota_X$ in Theorems \ref{thm6.5} and \ref{thm6.9}.
Moreover, the non-characteristic deformation and microlocal Morse lemmas apply to
these objects via their corresponding ordinary sheaves.

This article is organized as follows.
Section~\ref{sec2} recalls the definition of regularity for ind-sheaves.
Section~\ref{sec3} proves Proposition \ref{keyprop}, a key result used in Section \ref{sec4} to prove Theorem \ref{mainB}.
Theorem~\ref{mainA} is proved in Section~\ref{sec5}.
In Section~\ref{sec6}, we characterize constructible ind-sheaves under the functor $\iota$ (Theorems \ref{thm6.5} and \ref{thm6.9}), and relate them to regular holonomic $\D$-modules via the regular Riemann--Hilbert correspondence (Theorem \ref{thm6.10}). 
The appendix gives applications of Theorem \ref{thm6.5} to
a non-characteristic deformation lemma (Theorem \ref{mainD})
and a microlocal Morse lemma (Theorem \ref{mainE}) for regular ind-sheaves.

%\newpage

\section*{Acknowledgments}
The first author is very grateful to Professors Masatoshi Kitagawa and Taito Tauchi
for fruitful discussions and many useful comments.
The second author thanks Professors Naofumi Honda and Toshihiro Oshiba
for their helpful comments.

\section*{Funding}
The research of the first author was partially supported
by JSPS KAKENHI Grant Number 25K17237.
The research of the second author was partially supported
by JSPS KAKENHI Grant Number 26K24488, 
Japan Society for the Promotion of Science (JSPS).
This research was partially supported by the JST Moonshot R\&D Program (Grant Number JPMJMS2389).

\section*{Use of AI}
The second author used ChatGPT (Astra) to assist with improving the English presentation,
organizing and improving some proofs, and checking notation and consistency.
Although parts of the wording suggested by AI were adopted,
the final manuscript was prepared and typed by the authors, who take full responsibility for its contents.
All proofs were independently checked by the authors without the use of AI.

\section{Preliminaries on Ind-Sheaves}\label{sec2}
Throughout this paper, 
$\Bbbk$ is a field and $X$ is a real analytic or complex manifold.
We denote by $\I\Bbbk_X$ the category of ind-objects of $\Mod^c(\Bbbk_X)$ which are called ind-sheaves,
and by $\I_{\RRc}\Bbbk_X$ the category of ind-objects of $\Mod^c_{\RRc}(\Bbbk_X)$.
Here, $\Mod^c(\Bbbk_X)$ (resp.~$\Mod^c_{\RRc}(\Bbbk_X)$) is 
the category of sheaves (resp.~$\RR$-constructible sheaves) of $\Bbbk$-vector spaces on $X$
with compact support.
Note that there exists a natural exact embedding $\iota_X\colon\Mod(\Bbbk_X)\to\I\Bbbk_X$.
It has an exact left adjoint $\alpha_X\colon \I\Bbbk_X\to \Mod(\Bbbk_X)$, 
which has in turn an exact fully faithful
left adjoint functor $\beta_X\colon \Mod(\Bbbk_X)\to \I\Bbbk_X$.
We denote by $\BDC(\I\Bbbk_X)$ the derived category of $\I\Bbbk_X$,
and by $\BDC_{\I\RRc}(\I\Bbbk_X)$ the full subcategory of $\BDC(\I\Bbbk_X)$
consisting of objects whose cohomology sheaves belong to $\I_{\RRc}\Bbbk_X$.
%The derived categories of ind-sheaves admit the six operations
%\[
%\otimes,\quad \rihom,\quad \bfR f_\ast,\quad f^{-1},
%\quad \bfR f_{!!},\quad f^!,
%\]
%where $f\colon X\to Y$ is a morphism of manifolds.
%We also set $\rhom := \alpha_X\rihom$.

In this framework, Kashiwara--Schapira introduced the micro-support and regularity for ind-sheaves.
We recall the definition of regularity.
See \cite[Def.\:4.2, Prop.\:4.3]{KS03} for the definition and properties of the micro-support of ind-sheaves.
We use the functor defined in \cite[\S\:3]{KS03}:
$$J\colon \BDC(\I\Bbbk_X)\to {\rm Ind}(\BDC(\Mod^c(\Bbbk_X))),\
F\mapsto \Hom_{\BDC(\I\Bbbk_X)}(\cdot, F).$$
Then the notion of regularity for ind-sheaves is defined as follows.
\newpage

\begin{definition}[{\cite[Def.\:4.9]{KS03}}]\label{regularity}
We say that $F\in \BDC(\I\Bbbk_X)$ is regular along a locally closed conic subset $\Lambda\subset T^\ast X$ at $p\in T^\ast X$
if there exist a neighborhood $N$ of $\pi(p)$,
$F'\in \BDC(\I\Bbbk_X)$,
an open neighborhood $U$ of $p$,
a small filtrant category $I$,
and a functor $I\ni i\mapsto F_i\in \bfD^{[a,b]}(\Bbbk_X)$
such that
$F|_N\simeq F'|_N$,
$U\cap \Lambda$ is closed in $U$,
$\displaystyle J(F')\simeq``\varinjlim_{i\in I}"J(F_i)$
and $\SS(F_i)\cap U\subset \Lambda$ for any $i\in I$.

Moreover,we use similar terminology as in \cite[Def.\:4.9]{KS03}.
\end{definition}

%Then, we have the following properties of regularity for ind-sheaves.
%\begin{proposition}[{\cite[Prop.\:4.10]{KS03}}]\label{prop-regular}
%\begin{itemize}
%\item[{\rm (1)}]
%Let $F\in \BDC(\I\Bbbk_X)$.
%Then $F$ is regular along any locally closed conic subset $\Lambda\subset T^\ast X$ at each $p\notin \SS(F)$. 
%
%\item[{\rm (2)}\footnotemark]
%Let $F_1\to F_2\to F_3\xrightarrow{\ +1\ }$ be a distinguished triangle in $\BDC(\I\Bbbk_X)$.
%If $F_j$ and $F_k$ are regular along a locally closed conic subset $\Lambda$ at $p\in T^\ast X$,
%so is $F_i$ for $i,j,k\in\{1,2,3\}$ with $j\neq k$.
%
%\item[{\rm (3)}]
%For any $F\in \BDC(\Bbbk_X)$, $\iota_XF$ is regular.
%\end{itemize}
%\end{proposition}
%\footnotetext{Although the original statement \cite[Prop.\:4.10 (ii)]{KS03} is not formulated in a pointwise manner,
%the same proof also yields this pointwise version.}

The following results are well known to experts;
however, since we could not find its proof in the literature,
we prove it here.
This lemma means that if $F$ is regular along $\Lambda$ at $(x;0)\in T^\ast X$
then $F$ is so at $(x'; 0)$ for any $x'$ near $x$,
with the same data.

\begin{lemma}\label{lem-zero}
Let $F\in \BDC(\I\Bbbk_X)$ and $\Lambda\subset T^\ast X$ be a conic locally closed subset.
We assume that $F$ is regular along $\Lambda$ at $(x; 0)$,
and fix the data $(N, F',U, I, (F_i)_{i\in I})$ as in {\rm Definition \ref{regularity} (1)}.
Then there is a neighborhood $U_x$ of $x$ such that $F'|_{U_x}\simeq F|_{U_x}$
and for any $i\in I$ we have 
$$\SS(F_i)\cap T^\ast U_x\subset\Lambda.$$
\end{lemma}

\begin{proof}
By taking a local trivialization of $T^*X$ near $x$,
since $U$ is a neighborhood of $(x; 0)$,
there exist an open neighborhood $U_x\subset N$ of $x$ and $\varepsilon>0$ such that
$$\{(y; \xi)\in T^\ast U_x\ |\ \|\xi\| < \varepsilon\} \subset U.$$
Let $i\in I$ and $(y;\xi)\in \SS(F_i)\cap T^\ast U_x$.
If $\xi=0$, then $(y; 0)\in \SS(F_i)\cap U$, and hence $(y; 0)\in\Lambda$.
If $\xi\neq0$, we set $\displaystyle t := \frac{\varepsilon}{2\|\xi\|}$ and have $\| t\xi\| < \varepsilon$.
By the conicity of $\SS(F_i)$, we have 
$$(y; t\xi)\in \SS(F_i)\cap U\subset\Lambda.$$
Therefore, since $\Lambda$ is conic, we have $(y; \xi)\in\Lambda$.
The choices of $U_x$ and $\varepsilon$ are independent of $i$,
and we can take the same data as above.
\end{proof}

\begin{corollary}\label{cor2.3}
Let $F\in \BDC(\I\Bbbk_X)$ and $\Lambda\subset T^\ast X$ be a conic closed subset.
We assume that $F$ is regular along $\Lambda$ at $(x; 0)$,
and fix the data $(N, F',U, I, (F_i)_{i\in I})$ as in {\rm Definition \ref{regularity} }.
Then there exists an open neighborhood $V \subset N$ of $x$ such that
$\SS(F|_V)\subset\Lambda.$
\end{corollary}

\begin{proof}
By Lemma \ref{lem-zero}, shrinking $N$ if necessary,
we can fix the data $(N, F',U, I, (F_i)_{i\in I})$ as in {\rm Definition \ref{regularity} (1)}
such that $N$ is an open subset and $\SS(F_i)\cap T^\ast N\subset \Lambda$ for any $i\in I$.
Then we have
\[
\begin{aligned}
\SS(F|_N)
= \SS(F'|_N)
&= \SS(F')\cap T^\ast N\\
&\subset
 \left(\underset{i\in I}{\lim}\SS(F_i)\right)
 \cap T^\ast N\\ 
&\subset
 \overline{\bigcup_{i\in I}\SS(F_i)}
 \cap T^\ast N 
\subset
 \overline{\bigcup_{i\in I}
   \bigl(\SS(F_i)\cap T^\ast N\bigr)}
\subset \overline{\Lambda}
 = \Lambda,
\end{aligned}
\]
by \cite[Eq.\:(4.1)]{KS03}.
%Therefore, we have $\SS(F)\subset \Lambda$.
Taking $V$ to be this smaller neighborhood, we obtain $\SS(F|_V)\subset\Lambda\cap T^*V$.
\end{proof}

%\newpage

\section{Key Result of Theorem \ref{mainB}}\label{sec3}
In this section, we shall prove Proposition \ref{keyprop}, which is the key step in the proof of Theorem \ref{mainB}. 

Throughout this section, $\Lambda$ is a conic closed subanalytic isotropic subset of $T^\ast X$.
Let us denote by $a_x\colon \{x\}\hookrightarrow X$ the embedding for $x\in X$
and note that the functor $a_x^{-1}\colon \I\Bbbk_X\to \I\Bbbk_{\{x\}} = {\rm Ind}({\rm Vect}_\Bbbk)$
commutes with small filtrant ind-limits $``\varinjlim"$ and the functors $\iota, \alpha$.
See \cite[Props.\:4.3.1, 4.3.2]{KS01} for the details.
Note that $a_x^{-1}F$ differs from the ind-stalk ${}_xF\in \I\Bbbk_X$ of $F$ at $x$ defined
in \cite[p.\:70]{KS01}.

\begin{proposition}\label{keyprop}
Let $F\in \BDC(\I\Bbbk_X)$.
Assume that $F$ is regular along $\Lambda$ at $(x;0)\in T^\ast X$ with $x\in X$,
and there exists an open neighborhood $N$ of $x$ such that $F|_N\in \BDC_{\I\RRc}(\I\Bbbk_N)$
and $(\alpha_XF)|_N\simeq 0$.
Then there exists an open neighborhood $U_x$ of $x$ such that $F|_{U_x}\simeq 0$.
\end{proposition}

To prove this proposition, we prepare the following two lemmas.
We denote by ${\rm Vect}_{\Bbbk}^{\rm fd}$ the category of finite-dimensional $\Bbbk$-vector spaces.
\begin{lemma}\label{lem3.2}
Let $F\in {\rm Ind}({\rm Vect}_{\Bbbk}^{\rm fd})$.
Assume that $\alpha_{\{\rm pt\}}F\simeq 0$.
Then we have $F\simeq 0$.
\end{lemma} 

\begin{lemma}\label{lem3.3}
Let $F\in \BDC(\I\Bbbk_X)$.
Assume that $F$ is regular along $\Lambda$ at $(x;0)\in T^\ast X$ with $x\in X$
and there exists a neighborhood $N$ of $x$ such that $a_y^{-1}F\simeq 0$ for any $y\in N$.
Then there exists an open neighborhood $U_x$ of $x$ such that $F|_{U_x}\simeq 0$.
\end{lemma} 

Assuming these lemmas, we now prove Proposition \ref{keyprop}.
\begin{proof}[Proof of Proposition \ref{keyprop}]
Let $q\in \ZZ$ and $y\in N$ and we consider an ind-sheaf
$$H^q(a_y^{-1}F)\simeq (H^qF)_y\in {\rm Ind}({\rm Vect}_{\Bbbk}).$$
Since $F|_N\in \BDC_{\I\RRc}(\I\Bbbk_N)$, we have $H^qF|_N\in \I_{\RRc}\Bbbk_N$
and hence $(H^qF)_y\in {\rm Ind}({\rm Vect}_{\Bbbk}^{\rm fd})$.
By the assumption, we have $(\alpha_XF)_y\simeq 0$ and hence
$$\alpha_{\{y\}}(H^q(a_y^{-1}F)) \simeq \alpha_{\{y\}}((H^qF)_y) \simeq (H^q(\alpha_XF)_y)\simeq 0.$$
Thus, by Lemma \ref{lem3.2} we have $H^q(a_y^{-1}F)\simeq 0$ for any $q\in \ZZ$,
and hence $a_y^{-1}F\simeq 0$ for any $y\in N$.
Therefore, there exists an open neighborhood $U_x$ of $x$ such that $F|_{U_x}\simeq 0$
by Lemma \ref{lem3.3}.
\end{proof}

We shall now prove Lemma \ref{lem3.2}.
\begin{proof}[Proof of Lemma \ref{lem3.2}]
We fix a small filtrant inductive system $\{F_i\}_{i\in I}$ in ${\rm Vect}_{\Bbbk}^{\rm fd}$
such that $\displaystyle F = ``\varinjlim_{i\in I}"F_i$.
Then we have $\displaystyle \varinjlim_{i\in I}F_i \simeq 0$ by the assumption.
Let $i\in I$ and $S_i = \{{\bm x}^i_1,\ldots, {\bm x}^i_{N(i)}\}$ be a $\Bbbk$-basis of $F_i$.
Then for each $k\in \{1,\ldots, N(i)\}$, there exist an object $j_k$ and a morphism $i\to j_k$ in $I$
such that the image of ${\bm x}^i_k$ under the morphism $F_i\to F_{j_k}$ is equal to the zero vector ${\bm 0}$.
Since $I$ is filtrant, there exist an object $j$ and a morphism $i\to j$ in $I$
such that the morphism $i\to j$ factors through the morphism $i\to j_k$ for any $k\in \{1,\ldots, N(i)\}$.
Thus, the image of $S_i$ by the morphism $F_i\to F_j$ is equal to $\{\bm 0\}$,
and hence the morphism $F_i\to F_j$ is the zero map.
Therefore, we have $F\simeq 0$.
\end{proof}

We now prove Lemma \ref{lem3.3}.
To prove this, we prepare the following Lemma \ref{lem3.4}.
Let $K$ be a compact subanalytic neighborhood of $x$
and $\calT$ be a finite triangulation of $K$.
Then we denote by $\Mod_{\wc}^K(\Bbbk_X)$ the full abelian subcategory of weakly $\calT$-constructible sheaves
with support in $K$.
Note that the fully faithful functor $\Mod_{\wc}^K(\Bbbk_X)\hookrightarrow \Mod^c(\Bbbk_X)$ induces
an exact and fully faithful functor ${\rm Ind}(\Mod_{\wc}^K(\Bbbk_X))\hookrightarrow \I\Bbbk_X$
by \cite[Prop.\:1.1.3 (iii) and the discussion following Thm.\:1.3.1]{KS01}.
Moreover, we shall fix $x_\sigma\in \sigma$ for each $\sigma\in \calT$
and define $$\omega_\sigma\colon {\rm Ind}(\Mod_{\wc}^K(\Bbbk_X))\to {\rm Ind}({\rm Vect}_\Bbbk),\
``\varinjlim_{i\in I}"F_i\mapsto ``\varinjlim_{i\in I}"(F_i)_{x_\sigma}.$$

\begin{lemma}\label{lem3.4}
Let $F\in {\rm Ind}(\Mod_{\wc}^K(\Bbbk_X))$.
Assume that $\omega_\sigma(F)\simeq 0$ for any $\sigma\in \calT$.   
Then we have $F\simeq 0$.
Namely, $\{\omega_\sigma\}_{\sigma\in \calT}$ is jointly conservative.
\end{lemma}

\begin{proof}

We fix a small filtrant inductive system $\{F_i\}_{i\in I}$ in $\Mod_{\wc}^K(\Bbbk_X)$
such that $\displaystyle F = ``\varinjlim_{i\in I}"F_i$.
By the assumption, we have $\displaystyle ``\varinjlim_{i\in I}"(F_i)_{x_\sigma}\simeq 0$
for any $\sigma\in \calT$ .
Let $i\in I$.
Then for each $\sigma\in \calT$, there exist an object $j_\sigma$ and a morphism $i\to j_\sigma$ in $I$
such that the map $(F_i)_{x_\sigma}\to (F_{j_\sigma})_{x_\sigma}$ is the zero map.
Since $\calT$ is a finite set and $I$ is filtrant, 
there exist an object $j$ and a morphism $i\to j$ in $I$ such that
the morphism $i\to j$ factors through the morphism $i\to j_\sigma$ for any $\sigma\in \calT$,
and hence the map $(F_i)_{x_\sigma}\to (F_j)_{x_\sigma}$ is the zero map.
 Since each $F_i$ is weakly $\calT$-constructible,
 the morphism $\varphi_\sigma\colon F_i|_\sigma\to F_j|_\sigma$ is a morphism
 between locally constant sheaves.
Because the map $(\varphi_\sigma)_{x_\sigma}\colon (F_i)_{x_\sigma}\to (F_j)_{x_\sigma}$ is the zero map
and $\sigma$ is connected,
the morphism $\varphi_\sigma\colon F_i|_\sigma\to F_j|_\sigma$ is the zero morphism.
Therefore, since $\varphi_\sigma=0$ for any $\sigma\in \calT$ and 
the support of each $F_i$ is contained in $K$,
the morphism $F_i\to F_j$ is the zero morphism and hence we have $F\simeq 0$.
\end{proof}

\begin{proof}[Proof of Lemma \ref{lem3.3}]
By Lemma \ref{lem-zero}, shrinking $N$ if necessary,
we can fix the data $(N, F',U, I, (F_i)_{i\in I}, a, b)$ as in {\rm Definition \ref{regularity}}
such that $N$ is open and $\SS(F_i)\cap T^\ast N\subset \Lambda$ for any $i\in I$.

By \cite[Cor.\:8.3.22]{KS90},
there exists a locally finite subanalytic $\mu$-stratification $\calS$ of $N$
such that $$\Lambda \cap T^\ast N\subset \bigcup_{S\in\calS}T^\ast_SX.$$
Then $H^pF_i|_S$ is locally constant
for any $i\in I$, any $p\in \ZZ$ and any $S\in \calS$ by \cite[Prop.\:8.4.1]{KS90}.
Let us choose a compact subanalytic neighborhood $K$ of $x$ such that
$\operatorname{Int}K\subset K\Subset N$,
a finite triangulation $\calT$ of $K$ subordinate to the $\mu$-stratification $\calS$,
and $x_\sigma \in \sigma$ for each $\sigma\in \calT$.

We shall consider $F\otimes \Bbbk_K$.
Let $q\in \ZZ$.
By \cite[Thm.\:3.2 (iii), Prop.\:3.9]{KS03},
we have
$$
J(F\otimes \Bbbk_K) \simeq ``\varinjlim_{i\in I}"J(F_i\otimes \Bbbk_K),\hspace{17pt}
H^q(F\otimes \Bbbk_K) \simeq ``\varinjlim_{i\in I}"H^q(F_i)\otimes \Bbbk_K.
$$
Since $H^q(F_i) \otimes \Bbbk_K \in \Mod_{\wc}^K(\Bbbk_X)$,
we have $H^q(F\otimes \Bbbk_K)\in {\rm Ind}(\Mod_{\wc}^K(\Bbbk_X))$.
Moreover, for any $\sigma\in \calT$,
we have
$$
\omega_\sigma(H^q(F\otimes \Bbbk_K))
\simeq
``\varinjlim_{i\in I}"(H^q(F_i)\otimes \Bbbk_K)_{x_\sigma}
\simeq
``\varinjlim_{i\in I}"H^q(F_i)_{x_\sigma}
\simeq
H^q(F)_{x_\sigma}
\simeq
H^q(a_{x_\sigma}^{-1}F)
\simeq 0
$$
by the assumption.
Thus, we have $H^q(F\otimes \Bbbk_K)\simeq 0$ by Lemma \ref{lem3.4}.
Therefore,  we have $F\otimes \Bbbk_K\simeq 0$,
and hence it is enough to take $\Int\,K$ as $U_x$.
\end{proof}

By \cite[Prop.\:4.3.21]{KS01},
we have the following results as corollaries of Proposition \ref{keyprop} and Lemma \ref{lem3.3}.
\begin{corollary}
Let $F\in \BDC_{\I\RRc}(\I\Bbbk_X)$.
Assume that for any $x\in X$,
there exists a conic closed subanalytic isotropic subset $\Lambda_x\subset T^\ast X$
such that $F$ is regular along $\Lambda_x$ at $(x;0)\in T^\ast X$,
and there exists a neighborhood $U_x$ of $x$ such that $(\alpha_XF)|_{U_x}\simeq 0$.
Then we have $F\simeq 0$.
\end{corollary}

\begin{corollary}
Let $F\in \BDC_{\I\RRc}(\I\Bbbk_X)$.
Assume that for any $x\in X$,
there exists a conic closed subanalytic isotropic subset $\Lambda_x\subset T^\ast X$
such that $F$ is regular along $\Lambda_x$ at $(x;0)\in T^\ast X$ and $a_x^{-1}F\simeq 0$.
Then we have $F\simeq 0$.
\end{corollary}

We omit the proofs.

\section{Proof of Theorem \ref{mainB} and its Corollaries}\label{sec4}
In this section, we prove Theorem \ref{mainB} using Lemma \ref{lem-zero} and Proposition \ref{keyprop}.

\begin{proof}[Proof of Theorem \ref{mainB}]
By Lemma \ref{lem-zero}, shrinking $N$ if necessary,
we can fix the data $(N, F',U, I, (F_i)_{i\in I}, a, b)$ as in {\rm Definition \ref{regularity}}
such that $N$ is an open subset of $U_0$ and $\SS(F_i)\cap T^\ast N\subset \Lambda$ for any $i\in I$.
Then we have
$$\SS(F|_N)\subset \Lambda$$
by Corollary \ref{cor2.3}.
See \cite[Def.\:4.5]{KS03} for the definition of $lim$.
Moreover, there exists a locally finite subanalytic $\mu$-stratification $\calS$ of $N$
such that $$\Lambda \cap T^\ast N\subset \bigcup_{S\in\calS}T^\ast_SX$$
by \cite[Cor.\:8.3.22]{KS90}.
Since $\SS(\alpha_XF)\subset \SS(F)$ by \cite[Prop\:4.3 (iii)]{KS03},
we have
$$\SS((\alpha_XF)|_N) = \SS(\alpha_XF)\cap T^\ast N \subset \SS(F)\cap T^\ast N = \SS(F|_N)\subset \Lambda$$
and hence $(\alpha_XF)|_N$ is weakly $\RR$-constructible adapted to $\calS$ by \cite[Prop.\:8.4.1]{KS90} .
Moreover, $(\alpha_XF)|_N$ is $\RR$-constructible by the assumption.
Thus $(\alpha_XF)|_N\in\BDC_{\RRc, \Lambda}(\Bbbk_N)$.
 
We now consider a distinguished triangle
$$F\xrightarrow{u_F}\iota_X\alpha_XF\to C_F\xrightarrow{\ +1\ }$$ 
in $\BDC(\I\Bbbk_X)$.
Here, $C_F$ is the mapping cone of $u_F$.
Since $\alpha_X(u_F)$ is an isomorphism, we have
$\alpha_X C_F \simeq 0$. Moreover,
$C_F|_N \in D^b_{\I\mathbb{R}\text{-c}}(Ik_N)$, since both
$F|_N$ and $(\iota_X\alpha_X F)|_N$ belong to this category.
The latter object is regular along $\Lambda$ at $(x;0)\in T^\ast X$ because
$\SS((\alpha_X F)|_N) \subset \Lambda$.
Hence $C_F|_N$ is regular along $\Lambda$ at $(x;0)\in T^\ast X$
by %\cite[Prop.\:4.10 (ii)]{KS03}.
\cite[Prop.\:4.10 (ii)]{KS03}.
Therefore Proposition~\ref{keyprop} gives a neighborhood $U_x \subset N$
such that $C_F|_{U_x} \simeq 0$.
Consequently, $u_F|_{U_x}$ is an isomorphism.
% The proof is completed.
\end{proof}

\newpage
We now present several corollaries of Theorem \ref{mainB}.
For this purpose, we shall define the notion of finite regularity as follows.

\begin{definition}
We say that $F\in \BDC(\I\Bbbk_X)$ is finitely regular along $\Lambda$ at $x\in X$
if $F$ is regular along $\Lambda$ at $(x;0)$ and
$\dim_{\Bbbk} H^q(\alpha_XF)_x<\infty$
for any $q\in \ZZ$.
Moreover, we use similar terminology as in \cite[Def.\:4.9]{KS03}.
\end{definition}

\begin{corollary}\label{cor4.1}
Let $F\in \BDC(\I\Bbbk_X)$.
Assume that for any $x\in X$,
there exists an open neighborhood $U_x$ of $x$
such that $F|_{U_x}\in \BDC_{\I\RRc}(\I\Bbbk_{U_x})$,
there exists a conic closed subanalytic isotropic subset $\Lambda_x\subset T^\ast U_x$
such that $F$ is finitely regular along $\Lambda_x$ at $x$.
Then the canonical morphism $u_F\colon F\to \iota_X\alpha_XF$ is an isomorphism
and $\alpha_XF\in \BDC_{\RRc}(\Bbbk_X)$.
\end{corollary}

\begin{corollary}\label{cor4.2}
Let $F\in \BDC_{\I\RRc}(\I\Bbbk_X)$.
Assume that $\SS(F)$ is contained in a conic locally closed subanalytic isotropic subset,
and for any $x\in X$,
$F$ is finitely regular at $x$.
Then the canonical morphism $u_F\colon F\to \iota_X\alpha_XF$ is an isomorphism
and $\alpha_XF\in \BDC_{\RRc}(\Bbbk_X)$.
\end{corollary}

We omit the proofs of the above two corollaries.

\begin{corollary}
Let $F\in \BDC(\I\Bbbk_X)$ and $\Lambda\subset T^\ast X$ be a conic closed subanalytic isotropic subset,
and set $\Delta_F := {\rm Cone}(u_F\colon F\to \iota_X\alpha_XF)[-1]$.
Assume that $F\in \BDC_{\I\RRc}(\I\Bbbk_X)$, $\SS(F)\subset \Lambda$ and 
for any $x\in X$ and any $q\in\ZZ$, $\dim_{\Bbbk} H^q(\alpha_XF)_x<\infty$.
Then $F$ is regular along $\Lambda$ at $(x;0)\in T^\ast X$
if and only if
there exists an open neighborhood $U_x$ of $x$
such that $\Delta_F|_{U_x}\simeq 0$.
\end{corollary}

\begin{proof}
The ``only if\," part follows from Theorem \ref{mainB}.
We now prove the ``if\," part.
We assume that there exists an open neighborhood $U_x$ of $x$
such that $\Delta_F|_{U_x}\simeq 0$.
Then we have $u_F|_{U_x}\colon F|_{U_x}\simeq (\iota_X\alpha_XF)|_{U_x}$
and hence 
$$\SS(\alpha_XF)\cap T^\ast U_x
= \SS((\alpha_XF)|_{U_x})
= \SS((\iota_X\alpha_XF)|_{U_x})
=  \SS(F|_{U_x})
= \SS(F) \cap T^\ast U_x
%\subset \SS(F)
\subset \Lambda.$$
Therefore, we see that $F$ is regular along $\Lambda$ at $(x;0)\in T^\ast X$.
\end{proof}

We shall give applications of Corollary \ref{cor4.1} in Section \ref{sec6}.
At the end of this section, 
we introduce the following theorem, which is obtained by a similar argument to Theorem \ref{mainB}.

\begin{theorem}
Let $F\in \BDC(\I\Bbbk_X)$ and $\Lambda\subset T^\ast X$ be a conic closed subanalytic isotropic subset.
Assume that $F$ is regular along $\Lambda$ at $(x;0)\in T^\ast X$ with $x\in X$
and there exists an open neighborhood $U_0$ of $x$ such that
for any $y\in U_0$ and any $q\in \ZZ$ there exists a finite-dimensional $\Bbbk$-vector space $V_y^q$
such that $H^q(a_y^{-1}F)\simeq \iota_{\{y\}}V_y^q$.
Then there exists an open neighborhood $U_x$ of $x$ such that
$u_F|_{U_x}\colon F|_{U_x}\simto (\iota_X\alpha_XF)|_{U_x}$
and $(\alpha_XF)|_{U_x}\in \BDC_{\RRc, \Lambda}(\Bbbk_{U_x})$.
\end{theorem}

\begin{proof}
By the argument used in the proof of Theorem \ref{mainB}, shrinking $U_0$ if necessary, $(\alpha_XF)|_{U_0}$ is weakly $\RR$-constructible.
Moreover, $(\alpha_XF)|_{U_0}$ is $\RR$-constructible by the assumption.

Let $C_F$ be a mapping cone of the canonical morphism $u_F\colon F\to \iota_X\alpha_XF$.
% For any $y\in U_0$ and $q\in \ZZ$, we have
% $$H^q(a_y^{-1}\iota_X\alpha_XF)
% \simeq
% H^q(\iota_{\{y\}}\alpha_{\{y\}}a_y^{-1}F)
% \simeq
% H^q(\iota_{\{y\}}\alpha_{\{y\}}\iota_{\{y\}}V_y^q)
% \simeq
% H^q(\iota_{\{y\}}V_y^q)
% \simeq
% H^q(a_y^{-1}F)$$
% and hence $a_y^{-1}u_F$ is an isomorphism.
% Thus, we have $H^q(a_y^{-1}C_F)\simeq 0$ for any $y\in U_0$ and any $q\in \ZZ$.
% Therefore, there exists an open neighborhood $U_x$ of $x$ such that $C_F|_{U_x}\simeq 0$
% by Lemma \ref{lem3.3}
% and hence $u_F|_{U_x}\colon F|_{U_x}\simto (\iota_X\alpha_XF)|_{U_x}$.
The functor $a_y^{-1}$ commutes with $\alpha$ and $\iota$,
compatibly with the adjunction units. By exactness, the morphism
$H^q(a_y^{-1}u_F)$ identifies with
\[
u_{H^q(a_y^{-1}F)}:
H^q(a_y^{-1}F)
\longrightarrow
\iota_{\{y\}}\alpha_{\{y\}}H^q(a_y^{-1}F).
\]
This is an isomorphism because
$H^q(a_y^{-1}F)\simeq\iota_{\{y\}}V_y^q$ and the unit is an
isomorphism on the essential image of $\iota_{\{y\}}$.
Thus $a_y^{-1}u_F$ is an isomorphism, and hence
$a_y^{-1}C_F\simeq0$ for every $y\in U_0$.
%
%We now verify that 
As in the proof of Theorem \ref{mainB},
$C_F$ is regular along $\Lambda$ at $(x;0)\in T^\ast X$.
%By Lemma~2.3 and the fact that %$\SS(\alpha_XF)\subset \SS(F)$, after shrinking %$U_0$ if necessary, we have
%\[
%\SS((\alpha_XF)|_{U_0})
%\subset \SS(F|_{U_0})
%\subset \Lambda.
%\]
%Since $(\alpha_XF)|_{U_0}$ is an ordinary sheaf complex with
%micro-support contained in $\Lambda$, its image
%$(\iota_X\alpha_XF)|_{U_0}$ is regular along %$\Lambda$ at $(x;0)\in T^\ast X$.
%Together with the assumed regularity of $F$, the %distinguished triangle
%\[
%F \xrightarrow{u_F} \iota_X\alpha_XF
%\longrightarrow C_F \xrightarrow{\ +1\ }
%\]
%and Proposition~2.2(2) imply that $C_F$ is %regular along $\Lambda$
%at $(x;0)\in T^\ast X$.
%
Therefore, Lemma~3.3 gives an open neighborhood $U_x\subset U_0$
of $x$ such that $C_F|_{U_x}\simeq0$.
Consequently, $u_F|_{U_x}$ is an isomorphism.
\end{proof}

\section{Proof of Theorem \ref{mainA}}\label{sec5}
In this section, we prove Theorem \ref{mainA}.
Throughout this section, $X$ is a complex manifold and $\Bbbk = \CC$.
To prove Theorem \ref{mainA}, we prepare three lemmas.

We denote by $\BDCcoh(\D_X)$ the triangulated category of coherent $\D_X$-modules
and by $\Char(\M)$ the characteristic variety of $\M\in \BDCcoh(\D_X)$.

\begin{lemma}\label{lem5.1}
For any $\M\in \BDCcoh(\D_X)$,
we have $\SS(\Sol^{\rmt}_X(\M)) = \Char(\M)$.
\end{lemma}

\begin{proof}
$\SS(\Sol^\rmt_X(\M))\subset \Char(\M)$
follows from \cite[Thm.\:6.1]{KS03} by induction on the length of the complex.
Moreover, we have
$\Char(\M)=\SS(\Sol_X(\M))=\SS(\alpha_X\Sol^\rmt_X(\M))\subset \SS(\Sol^\rmt_X(\M))$.
Therefore, we have
$\SS(\Sol^{\rmt}_X(\M)) = \Char(\M)$.
\end{proof}

Since $\Char(\M)=\emptyset$ implies $\M\simeq 0$,
we have the following result.
\begin{corollary}\label{cor5.2}
Let $\M\in \BDCcoh(\D_X)$.
If $\Sol^{\rmt}_X(\M)\simeq 0$,
then we have $\M\simeq 0$.
\end{corollary}

\begin{lemma}\label{lem5.3}
For any $\M\in \BDCcoh(\D_X)$,
we have $\Sol_X^{\rmt}(\M)\in \BDC_{\I\RRc}(\I\CC_X)$.
\end{lemma}

\begin{proof}
By induction on the length of the complex, it suffices to treat the case where $\mathcal{M}$ is a coherent $\mathcal{D}_X$-module.
Moreover, since a coherent $\D_X$-module has locally free resolution of finite length and
the property of belonging to $\BDC_{\I\RRc}(\I\CC_X)$ is local property,
it is enough to show that the case where $\M=\D_X^{\oplus m}$ with $m\in \NN$.
In this case, it follows from $\Sol_X^\rmt(\D_X^{\oplus m}) \simeq (\SO_X^\rmt)^{\oplus m}$ and
the definition of $\SO_X^\rmt$.
\end{proof}

\begin{lemma}\label{lem5.4}
Let $\M\in \BDCcoh(\D_X)$.
Assume that there exists $F\in \BDC_{\CCc}(\CC_X)$ such that $\Sol^\rmt_X(\M)\simeq \iota_XF$.
Then we have $\M\in \BDCrh(\D_X)$.
\end{lemma}

\begin{proof}
For simplicity, we set
$$S := \Sol^\rmt \colon\bfD(\D_X)^\op\to \bfD(\I\CC_X),
\quad T := \rhom(\cdot, \SO_X^\rmt)\colon \bfD(\I\CC_X)\to \bfD(\D_X)^\op.$$
It is well known that $T\dashv S$.
We denote by
$$\varepsilon\colon {\rm id}_{\bfD(\I\CC_X)}\to S\circ T,\quad
\eta\colon T\circ S\to {\rm id}_{\bfD(\D_X)^\op}$$
the unit and counit morphisms. 
Then we have
$${\rm id}_{S(\M)}\simeq S(\eta_\M) \circ \varepsilon_{S(\M)} \colon 
S(\M) \to STS(\M)\to S(\M).$$

We fix $F\in \BDC_{\CCc}(\CC_X)$ such that $\Sol^\rmt_X(\M)\simeq \iota_XF$.
Then we have $$\M_{\reg} := \rhom(\iota_XF, \SO_X^\rmt)\in \BDCrh(\D_X)$$
and $$F\simeq \Sol_X(\M_\reg)$$
by the regular Riemann--Hilbert correspondence.
Moreover, we have $TS(\M)\simeq \M_\reg$ and
$$S(\M) = (\Sol^\rmt_X(\M))\simeq \iota_XF \simeq
\iota_X\Sol_X(\M_{\reg})\simeq \Sol_X^\rmt(\M_\reg) = S(\M_\reg).$$
Thus, it is enough to show that $\eta_\M\colon\M_\reg\to\M$ is an isomorphism.
Based on the above discussion, we have the following chain of equivalences:
\begin{align*}
\eta_\M \text{ is an isomorphism}
&\Longleftrightarrow
S(\eta_\M) \text{ is an isomorphism}\\
&\Longleftrightarrow
\varepsilon_{S(\M)} \text{ is an isomorphism}\\
&\Longleftrightarrow
\varepsilon_{S(\M_\reg)} \text{ is an isomorphism}\\
&\Longleftrightarrow
S(\eta_{\M_\reg}) \text{ is an isomorphism}\\
&\Longleftrightarrow
\eta_{\M_\reg} \text{ is an isomorphism},
\end{align*}
where the first and last equivalences follow from Corollary \ref{cor5.2}.
Since $\M_\reg\in \BDCrh(\D_X)$, $\eta_{\M_\reg}$ is an isomorphism.
Therefore, $\eta_\M$ is an isomorphism and hence we have $\M\in \BDCrh(\D_X)$.
\end{proof}

\begin{proof}[Proof of Theorem \ref{mainA}]
Since $\M$ is holonomic, $\alpha_X\operatorname{Sol}_X^t(M)\simeq\operatorname{Sol}_X(M)\in D^b_{\mathbb C\text{-c}}(\mathbb C_X)$.
(i) $\Rightarrow$ (ii) follows from \cite[Lem.\:7.4.11]{KS01} and \cite[Prop.\:4.10]{KS03}
as already seen in the Introduction.
(ii) $\Rightarrow$ (iii) is obvious.
(iii) $\Rightarrow$ (iv) follows from Lemmas \ref{lem5.1}, \ref{lem5.3} and Corollary \ref{cor4.2}.
(iv) $\Rightarrow$ (i) follows from Lemma \ref{lem5.4}.
\end{proof}

\section{Applications to Constructible Ind-Sheaves}\label{sec6}
In this section,
we characterize the essential image of constructible sheaves under the functor $\iota$ by using Corollary \ref{cor4.1}.

Let us recall the notion of $\RR$-constructible ind-sheaves.

\begin{definition}[{\cite[Def.\:4.7]{KS03}}]
\begin{itemize}
\item[(1)]
We denote by $\BDC_{\wRRc}(\I\Bbbk_X)$ the full triangulated subcategory of $\BDC_{\I\RRc}(\I\Bbbk_X)$ consisting of objects $F$ such that $\SS(F)$ is contained in a conic locally closed subanalytic isotropic subset of $T^\ast X$.

\item[(2)]
We denote by $\BDC_{\RRc}(\I\Bbbk_X)$ the full triangulated subcategory of $\BDC_{\wRRc}(\I\Bbbk_X)$ consisting of objects $F$ such that $\rhom(G,F)\in \BDC_{\RRc}(\Bbbk_X)$ for any $G\in \BDC_{\RRc}(\Bbbk_X)$.
\end{itemize}
\end{definition}
Note that functors $\iota_X$ and $\alpha_X$ induce functors:
\[\xymatrix@M=5pt@C=35pt{
\BDC_{\wRRc}(\Bbbk_X)
& \BDC_{\wRRc}(\I\Bbbk_X), \ar@{->}[l]_-{\alpha_X}}\hspace{17pt}
\xymatrix@M=5pt@C=55pt{
\BDC_{\RRc}(\Bbbk_X)\ar@<1.0ex>@{->}[r]^-{\iota_X}
& \BDC_{\RRc}(\I\Bbbk_X). \ar@<1.0ex>@{->}[l]^-{\alpha_X}}\]
However, the functor $\iota_X$ does not induce a functor from $\BDC_{\wRRc}(\Bbbk_X)$ to $\BDC_{\wRRc}(\I\Bbbk_X)$.

Let $\Lambda\subset T^\ast X$ be a conic closed subanalytic isotropic subset.
Then we have:
\begin{proposition}
Let $F\in \BDC(\I\Bbbk_X)$.
Assume that $F$ is weakly $\RR$-constructible and finitely regular along $\Lambda$ at any $x\in X$.
Then, $F$ is $\RR$-constructible and $\SS(F)\subset \Lambda$. 
\end{proposition}

\begin{proof}
Since $\alpha_XF\in \BDC_{\wRRc}(\Bbbk_X)$,
we have $\alpha_XF\in \BDC_{\RRc}(\Bbbk_X)$
by the assumption of finiteness of stalks of cohomologies.
For any $G\in \BDC_{\RRc}(\Bbbk_X)$,
we have
$$\rhom(G,F)
\simto
\rhom(\iota_XG,\iota_X\alpha_XF)
\simeq \rhom(G,\alpha_XF)\in\BDC_{\RRc}(\Bbbk_X),$$
where we used Corollary \ref{cor4.1} (resp.~\cite[Prop.\:8.4.10 (ii)]{KS90}) in the first (resp.~last) step.
Therefore, $F$ is $\RR$-constructible.

The second assertion follows from Corollary \ref{cor2.3}.
\end{proof}

We denote by $\BDC_{\frRRc, \Lambda}(\I\Bbbk_X)$ the full subcategory of $\BDC_{\RRc}(\I\Bbbk_X)$
consisting of objects $F$ which are finitely regular along $\Lambda$ at each $x\in X$.
Then we have the following result by
Corollary \ref{cor4.1} and the fact that $\alpha_X\circ \iota_X\simeq {\rm id}$ (see \cite[\S\:4]{KS01} for the details).

\begin{theorem}\label{thm6.3}
The functors $\iota_X$ and $\alpha_X$ induce an equivalence of categories
\[\xymatrix@M=5pt@C=55pt{
\BDC_{\RRc, \Lambda}(\Bbbk_X)\ar@<1.0ex>@{->}[r]^-{\iota_X}\ar@{}[r]|-\sim
& \BDC_{\frRRc, \Lambda}(\I\Bbbk_X). \ar@<1.0ex>@{->}[l]^-{\alpha_X}
}\]
\end{theorem}

By Corollary \ref{cor4.1}, we have:
\begin{proposition}
For any $F\in \BDC_{\wRRc}(\I\Bbbk_X)$,
$F$ is finitely regular at any $x\in X$
if and only if
for any $x\in X$ there exist an open neighborhood $U_x$ of $x$ and
a conic closed subanalytic isotropic subset $\Lambda_x\subset T^\ast U_x$
such that $F$ is finitely regular along $\Lambda_x$ at $x$.
\end{proposition}

We also define $\BDC_{\frRRc}(\I\Bbbk_X)$ as the full subcategory of $\BDC_{\RRc}(\I\Bbbk_X)$ consisting of objects $F$ which are finitely regular at each $x\in X$.
Then the essential image of $\BDC_{\RRc}(\Bbbk_X)$ under the functor $\iota_X$ is nothing but this category.
Namely, we obtain the following equivalence of categories as a consequence of Theorem \ref{thm6.3}:
\begin{theorem}\label{thm6.5}
The functors $\iota_X$ and $\alpha_X$ induce an equivalence of categories
\[\xymatrix@M=5pt@C=41pt{
\BDC_{\RRc}(\Bbbk_X)\ar@<1.0ex>@{->}[r]^-{\iota_X}\ar@{}[r]|-\sim
& \BDC_{\frRRc}(\I\Bbbk_X).\ar@<1.0ex>@{->}[l]^-{\alpha_X} }\]
\end{theorem}

For the rest of this section, let $X$ be a complex manifold, $\Bbbk = \CC$ and 
Let $\Lambda\subset T^\ast X$ be a $\CC^\times$-conic closed subanalytic isotropic subset.
We now introduce the notion of $\CC$-constructibility for ind-sheaves.

\begin{definition}
\begin{itemize}
\item[(1)]
We denote by $\BDC_{\wCCc}(\I\Bbbk_X)$ the full triangulated subcategory of $\BDC_{\wRRc}(\I\Bbbk_X)$ consisting of objects $F$ whose micro-support is contained in a $\CC^\times$-conic closed subanalytic isotropic subset of $T^\ast X$.

\item[(2)]
We set $\BDC_{\CCc}(\I\Bbbk_X) := \BDC_{\wCCc}(\I\Bbbk_X)\cap \BDC_{\RRc}(\I\Bbbk_X).$
\end{itemize}
\end{definition}

Although, the functor $\iota_X$ does not induce a functor from $\BDC_{\wCCc}(\Bbbk_X)$ to $\BDC_{\wCCc}(\I\Bbbk_X)$, we have:
\begin{proposition}\label{prop6.7}
The functors $\iota_X$ and $\alpha_X$ induce functors:
\[\xymatrix@M=5pt@C=35pt{
\BDC_{\wCCc}(\Bbbk_X)
& \BDC_{\wCCc}(\I\Bbbk_X), \ar@{->}[l]_-{\alpha_X}}\hspace{17pt}
\xymatrix@M=5pt@C=55pt{
\BDC_{\CCc}(\Bbbk_X)\ar@<1.0ex>@{->}[r]^-{\iota_X}
& \BDC_{\CCc}(\I\Bbbk_X). \ar@<1.0ex>@{->}[l]^-{\alpha_X}}\]
\end{proposition}

\begin{proof}
Let $F\in \BDC_{\CCc}(\Bbbk_X)$.
Then $F\in \BDC_{\RRc}(\Bbbk_X)$ and $\SS(F)$ is $\CC^\times$-conic by \cite[Def.\:8.5.6 (i)]{KS90},
and hence $\iota_XF\in \BDC_{\RRc}(\I\Bbbk_X)$ and $\SS(\iota_XF) = \SS(F)$ is $\CC^\times$-conic.
Thus, we have $\iota_XF\in \BDC_{\CCc}(\I\Bbbk_X)$.

Let $F\in \BDC_{\wCCc}(\I\Bbbk_X)$.
Then $\alpha_XF\in \BDC_{\wRRc}(\Bbbk_X)$ and $\SS(\alpha_XF)\subset \SS(F)$ is contained in a $\CC^\times$-conic closed subanalytic isotropic subset of $T^\ast X$,
and hence $\SS(\alpha_XF)$ is $\CC^\times$-conic by \cite[Thm.\:8.5.5]{KS90}.
Thus, we have $\alpha_XF\in \BDC_{\wCCc}(\Bbbk_X)$.

The functor $\alpha_X$ preserves $\CC$-constructibility since
$\BDC_{\CCc}(\I\Bbbk_X) := \BDC_{\wCCc}(\I\Bbbk_X)\cap \BDC_{\RRc}(\I\Bbbk_X)$ and
it preserves $\RR$-constructibility.
\end{proof}

\newpage
We set $\BDC_{\frCCc, \Lambda}(\I\CC_X) := \BDC_{\frRRc, \Lambda}(\I\CC_X)\cap \BDC_{\CCc}(\I\CC_X)$. 
Then, we obtain the following result from Theorem \ref{thm6.3} and Proposition \ref{prop6.7}.
\begin{theorem}\label{thm6.8}
The functors $\iota_X$ and $\alpha_X$ induce an equivalence of categories
\[\xymatrix@M=5pt@C=55pt{
\BDC_{\CCc, \Lambda}(\CC_X)\ar@<1.0ex>@{->}[r]^-{\iota_X}\ar@{}[r]|-\sim
& \BDC_{\frCCc, \Lambda}(\I\CC_X). \ar@<1.0ex>@{->}[l]^-{\alpha_X}
}\]
\end{theorem}

We also define $\BDC_{\frCCc}(\I\CC_X)$ as the full subcategory of $\BDC_{\CCc}(\I\CC_X)$ consisting of objects $F$ which are finitely regular at each $x\in X$.
Note that $$\BDC_{\frCCc}(\I\CC_X) = \BDC_{\frRRc}(\I\CC_X)\cap \BDC_{\CCc}(\I\CC_X).$$
Then the essential image of $\BDC_{\CCc}(\CC_X)$ under the functor $\iota_X$ is nothing but this category.
Namely, we obtain the following equivalence of categories as a consequence of Theorem \ref{thm6.8}:
\begin{theorem}\label{thm6.9}
The functors $\iota_X$ and $\alpha_X$ induce an equivalence of categories
\[\xymatrix@M=5pt@C=41pt{
\BDC_{\CCc}(\CC_X)\ar@<1.0ex>@{->}[r]^-{\iota_X}\ar@{}[r]|-\sim
& \BDC_{\frCCc}(\I\CC_X).\ar@<1.0ex>@{->}[l]^-{\alpha_X} }\]
\end{theorem}

Moreover, we have the following equivalences of categories by Lemma \ref{lem5.1}, Theorems \ref{thm6.8} and \ref{thm6.9} and the regular Riemann--Hilbert correspondence.
We denote by $\BDC_{{\rm rh}, \Lambda}(\D_X)$ the full subcategory of $\BDCrh(\D_X)$
consisting of objects whose characteristic varieties are contained in $\Lambda$.
\begin{theorem}\label{thm6.10}
There exist an equivalence of categories:
\[\xymatrix@M=7pt@C=61pt{
\BDC_{{\rm rh}, \Lambda}(\D_X)^{\mathrm{op}}\ar@<1.0ex>@{->}[r]^-{\Sol^{\rmt}_X}\ar@{}[r]|-\sim
& \BDC_{\frCCc, \Lambda}(\I\CC_X),\ar@<1.0ex>@{->}[l]^-{\rhom(\cdot, \SO_X^\rmt)}}\]
\[
\xymatrix@M=7pt@C=41pt{
\BDCrh(\D_X)^{\mathrm{op}}
\ar@<1.0ex>@{->}[r]^-{\Sol^{\rmt}_X}
\ar@{}[r]|-\sim
&
\BDC_{\frCCc}(\I\CC_X).
\ar@<1.0ex>@{->}[l]^-{\rhom(\cdot,\SO_X^\rmt)}
}\]
\end{theorem}

Furthermore, these equivalences of triangulated categories induce the following exact equivalences of Abelian categories:
\[
\xymatrix@M=7pt@C=61pt{
\Mod_{{\rm rh},\Lambda}(\D_X)^{\mathrm{op}}
\ar@<1.0ex>@{->}[r]^-{\Sol^{\rmt}_X(\cdot)[d_X]}
\ar@{}[r]|-\sim
&
\{F\in\BDC_{\frCCc, \Lambda}(\I\CC_X)\mid\alpha_XF\in\Perv(\CC_X)\}
\ar@<1.0ex>@{->}[l]^-{\rhom(\cdot,\SO_X^\rmt)[d_X]}
}
\]
\[
\xymatrix@M=7pt@C=41pt{
\Modrh(\D_X)^{\mathrm{op}}
\ar@<1.0ex>@{->}[r]^-{\Sol^{\rmt}_X(\cdot)[d_X]}
\ar@{}[r]|-\sim
&
\{F\in\BDC_{\frCCc}(\I\CC_X)
  \mid\alpha_XF\in\Perv(\CC_X)\}
\ar@<1.0ex>@{->}[l]^-{\rhom(\cdot,\SO_X^\rmt)[d_X]}
}
\]

\section{Appendix}
In this section,
we present a non-characteristic deformation lemma and a microlocal Morse lemma
for regular ind-sheaves as applications of Theorem \ref{thm6.5}.

The following result is a non-characteristic deformation lemma for regular ind-sheaves.
\begin{theorem}\label{mainD}
Let $X$ be a real analytic manifold,
$F\in \BDC_{\frRRc}(\I\Bbbk_X)$ and
$\{U_t\}_{t\in \RR}$ a family of open subsets of $X$.
We assume the following conditions\footnote{For the non-characteristic
deformation lemma, the corrected definition is
$Z_s=\bigcap_{t>s}\overline{U_t\setminus U_s}$;
see the correction to \cite[Proposition~2.7.2]{KS90} in
%\url{https://webusers.imj-prg.fr/~pierre.schapira/BooksMono/Errata.pdf}
https://webusers.imj-prg.fr/~pierre.schapira/BooksMono/Errata.pdf.}:
\begin{itemize}
\item[\rm (i)]
For any $t\in \RR$, $\displaystyle U_t = \bigcup_{s<t}U_s$,

\item[\rm (ii)]
For any $s,t\in \RR$ with $s\leq t$, $\var{U_t\setminus U_s}\cap \supp(F)$ is compact,

\item[\rm (iii)]
For any $s,t\in \RR$ with $s\leq t$ and $x\in Z_s\setminus U_t$,
$(\bfR\Gamma_{X\setminus U_t}F)|_{\{x\}}\simeq 0$.
Here, $\displaystyle Z_s := \bigcap_{t>s}\var{U_t\setminus U_s}$.
\end{itemize}

Then, for any $t\in \RR$, the morphism
$$\rHom(\Bbbk_{\cup_{s\in \RR}U_s}, F)\to \rHom(\Bbbk_{U_t}, F)$$
induced by the restriction morphism is an isomorphism.
\end{theorem}

\begin{proof}
Since $\supp(\alpha_XF)\subset \supp(F)$,
for any $s,t\in \RR$ with $s\leq t$, $\var{U_t\setminus U_s}\cap \supp(\alpha_XF)$ is compact. 

Let $s,t\in \RR$ with $s\leq t$ and $x\in Z_s\setminus U_t$.
Since $u_F\colon F\simto \iota_X\alpha_XF$,
we have
\begin{align*}
\iota_{\{x\}}((\bfR\Gamma_{X\setminus U_t}\alpha_XF)|_{\{x\}})
&\simeq
(\iota_X\bfR\Gamma_{X\setminus U_t}\alpha_XF)|_{\{x\}}\\
&\simeq
(\bfR\Gamma_{X\setminus U_t}\iota_X\alpha_XF)|_{\{x\}}\\
&\simot
(\bfR\Gamma_{X\setminus U_t}F)|_{\{x\}}
\simeq 0
\end{align*}
and hence
$$(\bfR\Gamma_{X\setminus U_t}\alpha_XF)_{\{x\}}\simeq 0.$$
Thus, we can apply the non-characteristic deformation lemma \cite[Prop.\:2.7.2]{KS90}
to the sheaf $\alpha_XF$
and hence for any $t\in \RR$,
the restriction map
$$\bfR\Gamma\left(\bigcup_{s}U_s; \alpha_XF\right)\to \bfR\Gamma(U_t; \alpha_XF)$$
is an isomorphism.
In general, for any open subset $U\subset X$,
we have
\begin{align*}
\bfR\Gamma(U; \alpha_XF)
\simeq \rHom(\Bbbk_U, \alpha_XF)
\simeq \rHom(\iota_X\Bbbk_U, \iota_X\alpha_XF)
\simeq \rHom(\Bbbk_U, F).
\end{align*}
Therefore, the morphism
$$\rHom(\Bbbk_{\cup_{s\in \RR}U_s}, F)\to \rHom(\Bbbk_{U_t}, F)$$
induced by the restriction morphism is an isomorphism.
\end{proof}

We recall that ${}_Z\Bbbk := \beta_X\Bbbk_Z$ for a locally closed subset $Z\subset X$. 
See \cite[Prop.\:4.2.14 (i)]{KS01} for the details.
Then the following result is a microlocal Morse lemma for regular ind-sheaves.
\begin{theorem}\label{mainE}
Let $X$ be a real analytic manifold,
$a\in\RR, b\in \RR\cup\{\infty\}$ with $a<b$
and $\varphi\colon X\to \RR$ a $\mathcal{C}^1$-function.
We assume that $F\in \BDC(\I\Bbbk_X)$ satisfies the following conditions:
\begin{itemize}
\item[\rm (i)]
$\varphi$ is proper on $\supp(F)$,

\item[\rm (ii)]
For any $x\in \varphi^{-1}([a,b))$, $(x, d\varphi(x))\notin \SS(F)$.
\end{itemize}
Then the morphism
$$\rHom({}_{\varphi^{-1}((-\infty, b))}\Bbbk, F)\to \rHom({}_{\varphi^{-1}((-\infty, a))}\Bbbk, F)$$
induced by the restriction morphism is an isomorphism.
Moreover, if $F\in \BDC_{\frRRc}(\I\Bbbk_X)$ then 
the morphism
$$\rHom(\Bbbk_{\varphi^{-1}((-\infty, b))}, F)\to \rHom(\Bbbk_{\varphi^{-1}((-\infty, a))}, F)$$
induced by the restriction morphism is an isomorphism.
\end{theorem}

\begin{proof}
Since $\supp(\alpha_XF)\subset \supp(F)$ and $\SS(\alpha_XF)\subset \SS(F)$,
$\varphi$ is proper on $\supp(\alpha_XF)$ and
for any $x\in \varphi^{-1}([a,b))$ we have $(x;d\varphi(x))\notin \SS(\alpha_XF)$.
Thus, we can apply the microlocal Morse lemma \cite[Cor.\:5.4.19.]{KS90} to the sheaf $\alpha_XF$
and hence the restriction map
$$\bfR\Gamma_{\varphi^{-1}((-\infty, b))}(X; \alpha_XF)\to
\bfR\Gamma_{\varphi^{-1}((-\infty, a))}(X; \alpha_XF)$$
is an isomorphism.
In general, for any open subset $U\subset X$,
we have
$$\bfR\Gamma_{U}(X; \alpha_XF)
\simeq
\bfR\Gamma(U; \alpha_XF)
\simeq
\rHom(\Bbbk_U, \alpha_XF)
\simeq
\rHom(\beta_X\Bbbk_U, F)
\simeq
\rHom({}_U\Bbbk, F).$$
Therefore, the morphism
$$\rHom({}_{\varphi^{-1}((-\infty, b))}\Bbbk, F)\to \rHom({}_{\varphi^{-1}((-\infty, a))}\Bbbk, F)$$
induced by the restriction morphism is an isomorphism.

Moreover, we assume that $F\in \BDC_{\frRRc}(\I\Bbbk_X)$.
Then we have $u_F\colon F\simto \iota_X\alpha_XF$
and hence
$$\rHom({}_U\Bbbk, F)
\simeq \rHom(\Bbbk_U, \alpha_XF)
\simeq \rHom(\iota_X\Bbbk_U, \iota_X\alpha_XF)
\simot \rHom(\iota_X\Bbbk_U, F).$$
Therefore, the morphism
$$\rHom(\Bbbk_{\varphi^{-1}((-\infty, b))}, F)\to \rHom(\Bbbk_{\varphi^{-1}((-\infty, a))}, F)$$
is an isomorphism.
\end{proof}

% \begin{thebibliography}{99}


\begin{thebibliography}{99}
\bibitem[Bj{\"o}93]{Bjo93}
Jan-Erik Bj{\"o}rk,
Analytic $D$-modules and applications, volume 247 of {
  Mathematics and its Applications},
Kluwer Academic Publishers Group, Dordrecht, 1993.

\bibitem[HTT]{HTT}
Ryoshi Hotta, Kiyoshi Takeuchi, and Toshiyuki Tanisaki,
{{$D$}-modules, perverse sheaves, and representation theory},
  volume 236 of {Progress in Mathematics},
  Birkh{\"a}user Boston, 2008.

%\bibitem[Ka95]{Ka95}
%Masaki Kashiwara, Algebraic study of systems of partial differential equations, M\UTF{00E9}m. Soc. Math. France 63, 1995.

\bibitem[KS90]{KS90}
Masaki Kashiwara, Pierre Schapira,
Sheaves on manifolds, Grundlehren der Mathematischen Wissenschaften,
{\bf292}, Springer-Verlag, 1990.

\bibitem[KS01]{KS01}
Masaki Kashiwara, Pierre Schapira,
Ind-sheaves,
Ast\'{e}risque, {\bf 271}, 2001.

\bibitem[KS03]{KS03}
Masaki Kashiwara, Pierre Schapira, Microlocal study of ind-sheaves I: micro-support and regularity,
Ast\'{e}risque, {\bf 284}, 143--164, 2003.

\bibitem[KS16]{KS16}
Masaki Kashiwara, Pierre Schapira,
Regular and irregular holonomic $\D$-modules, 
 London Mathematical Society Lecture Note Series, {\bf433},
 Cambridge University Press, 2016.

\bibitem[Mar]{Mar}
Ana Rita Martins,
Temperate holomorphic solutions and regularity of holonomic $\D$-modules on curves,
Rend. Semin. Mat. Univ. Padova, \textbf{118}, 35--48, 2007.
\end{thebibliography}
\end{document}